\documentclass[10pt]{article}
\usepackage{amsmath}
\usepackage{latexsym}
\usepackage{amsfonts}
\usepackage{latexsym}

\def\Bbb{\mathbb}
\def\bcp{\mathbb C \mathbb P}
\def\CC{\mathbb C}

\def\dim{\mbox{dim}}
\def\Tor{\mbox{torsion}}

\newtheorem{main}{Theorem}

\newtheorem{defn}{Definition}
\newtheorem{thm}{Theorem}
\newtheorem{prop}[thm]{Proposition}

\newtheorem{lem}[thm]{Lemma}

\def\ro{\stackrel{\circ}{r}}

\newenvironment{proof}{\medskip \noindent
{\bf Proof.}}{\hfill \rule{.5em}{1em}
\\}

\newenvironment{rmk}{\mbox{ }\\{\bf  Remark}\mbox{ }}{
\hfill $\Box$\mbox{}\bigskip}
\def\ZZ{{\mathbb Z}}
\def\RR{{\mathbb R}}

\begin{document}
\def\NN{\mathbb N}
\def\ZZ{\mathbb Z}
\def\QQ{\mathbb Q}
\def\RR{\mathbb R}
\def\CC{\mathbb C}
\def\SS{\mathbb S}
\def\PP{\mathbb P}
\def\VV{\mathbb V}
\def\ro{\mathring{r}}
\def\ve{\varepsilon}
\def\bcp{\mathbb C \mathbb P}
\def\cpb{\overline{\mathbb C \mathbb P^2}}
\def\KK{\cal K}
\def\EE{\cal E}
\def\LL{\cal L}
\def\OO{\cal O} 
\def\XX{\cal X}
\def\pp{\cal P}
\def\mm{\cal M}
\def\cc{\cal C}
\def\yy{\cal Y}
\def\zz{\cal Z}
\def\dd{\Delta}
\def\ss{\Sigma}
\def\T{\Theta}

\title{On Normalized Ricci Flow and Smooth Structures on Four-Manifolds with $b^+=1$}

\author{Masashi Ishida, Rare{\c s} R{\u a}sdeaconu, and Ioana {\c S}uvaina}

\date{}

\maketitle

\begin{abstract}
We find an obstruction to the existence of non-singular solutions to the normalized Ricci flow on four-manifolds with $b^+=1$. By using this obstruction, we study the relationship between the existence or non-existence of non-singular solutions of the normalized Ricci flow and exotic smooth structures on the topological 4-manifolds ${\mathbb C}{P}^2 \# {\ell} \overline{{\mathbb C}{P}^2}$, where $5 \leq \ell \leq 8$. 
\end{abstract}


\section{Introduction}

Let $X$ be a closed oriented Riemannian manifold of dimension $n\geq 3.$
The normalized Ricci flow on  $X$  is the following evolution equation:
\begin{eqnarray*}\label{Ricci}
 \frac{\partial }{\partial t}{g}=-2{Ric}_{g} + \frac{2}{n}\Big(\frac{{\int}_{X} {s}_{g} d{\mu}_{g}}{{\int}_{X}d{\mu}_{g}} \Big) {g}, 
\end{eqnarray*}
where ${Ric}_{g}, s_g$ are the Ricci curvature and the scalar curvature of the evolving Riemannian metric $g$ and $d{\mu}_{g}$ is the volume measure with respect to $g$. Recall that a one-parameter family of metrics $\{g(t)\}$, where $t \in [0, T)$ for some $0<T\leq \infty$, is called a solution to the normalized Ricci flow if it satisfies the above equation at all $x \in X$ and $t \in [0, T)$. A solution $\{g(t)\}$ on a time interval $[0, T)$ is said to be maximal if it cannot be extended past time $T$. In this paper we are interested in solutions which are particularly nice. The following definition was first introduced and studied by Hamilton \cite{ha-0, c-c}: 
\begin{defn}\label{non-sin}
A maximal solution $\{g(t)\}$, $t \in [0, T)$, to the normalized Ricci flow on $X$ is called non-singular if $T=\infty$ and the Riemannian curvature tensor $Rm_{g(t)}$ of $g(t)$ satisfies $\sup_{X \times [0, \infty)}|Rm_{g(t)}| < \infty$.
\end{defn} 

Fang and his collaborators \cite{fz-1} pointed out that for a $4-$manifold with
negative Perelman invariant \cite{p-1, lott}, which is equivalent to the Yamabe
invariant in this situation \cite{A-ishi-leb-3}, the existence of non-singular
solutions forces a topological constraint on the $4-$manifold. On the other hand, the first author proved \cite{ism} that, in dimension four, the existence of the non-singular solutions is fundamentally related to the smooth structure considered. An important ingredient in his theorems was the non-triviality of the Seiberg-Witten invariant. In the case when the underlying manifold has $b^+\geq2,$ this is a diffeomorphism invariant. However, when $b^+=1$ the invariant depends on a choice of an orientation on $H^2(X,\ZZ)$ and $H^1(X,\RR).$ The obstructions in \cite{ism,i-s} are for manifolds with $b^+\geq2.$ We extend these results to the case $b^+=1,$ and we have the following:
\begin{main}\label{main-A}
 Let $X$ be a closed oriented smooth 4-manifold with  ${b}^+(X) = 1$ and $2 \chi(X) + 3\tau(X) > 0$. Assume that $X$ has a non-trivial Seiberg-Witten invariant.  Then, there do not exist non-singular solutions to the normalized Ricci flow on $X \# k \overline{ {\mathbb C}{P}^{2} }$ if 
\begin{eqnarray*}
k \geq \frac{1}{3} \Big( 2\chi(X) + 3 \tau(X) \Big). 
\end{eqnarray*}
\end{main}

We prove this result in a slightly more general setting, see Theorem \ref{main-C} in Section \ref{sec-4} below. Using Theorem \ref{main-A}, we study manifolds with small topology and emphasize how the change of smooth structure reflects on the existence or non-existence of solutions of the normalized Ricci flow:
\begin{main} \label{main-B}
For $5 \leq \ell \leq 8$, the topological 4-manifold $M:={\mathbb C}{P}^2 \# {\ell} \overline{{\mathbb C}{P}^2}$ satisfies the following properties:
\begin{itemize}
\item [ 1. ]  $M$ admits a smooth structure of positive Yamabe invariant on which there exists a non-singular solution to the normalized Ricci flow.
\item [ 2. ] $M$ admits a smooth structure of negative Yamabe invariant on which there exist  non-singular solutions to the normalized Ricci flow.
\item [ 3. ] $M$ admits infinitely many distinct smooth structures all of which have negative Yamabe invariant and on which there are no non-singular solutions to the normalized Ricci flow for any initial metric.
\end{itemize}  
\end{main}

\bigskip
\noindent
{\bf Acknowledgments.} 
We would like to express our deep gratitude to Claude LeBrun for his warm encouragement. The first author is partially supported by the Grant-in-Aid for Scientific Research (C), Japan Society for the Promotion of Science, No. 20540090. The second author would like to thank CNRS and IRMA Strasbourg for support and the excellent conditions provided while this work was completed.

\section{Polarized 4-manifolds and Seiberg-Witten invariants}

Let $X$ be a closed oriented smooth 4-manifold. Any Riemannian metric $g$
 on $X$ gives rise to a decomposition $H^2(X, {\mathbb R})=
{\mathcal H}^+_{g} \oplus {\mathcal H}^-_{g}$, where 
${\mathcal H}^+_{g},$ ${\mathcal H}^-_{g}$  consists 
of cohomology classes for which the harmonic representative is 
$g$-self dual or $g$-anti-self-dual, respectively. Notice that 
$b^+(X):=\dim {\mathcal H}^+_{g}$ is a non-negative integer 
which is independent of the metric $g$. In this article, we 
always assume that $b^+(X) \geq 1$ and mainly consider the case 
when $b^+(X)=1$. For a fixed $b^+(X)$-dimensional subspace 
$H \subset H^2(M, {\mathbb R})$ on which the intersection form is positively defined,
 we consider the set of all Riemannian 
metrics $g$ for which ${\mathcal H}^+_{g} = H$ is satisfied. 
The Riemannian metric $g$ satisfying this property is called a
$H$-adapted metric. Under the assumption that there is at least 
one $H$-adapted metric, $H$ is called a polarization of $X$ and 
we call $(X, H)$ a polarized 4-manifold following \cite{leb-2}. 
For any given element $\alpha \in H^2(X, {\mathbb R})$ and a 
polarization $H$ of $X$, we use $\alpha^{+}$ to 
denote the orthogonal projection of $\alpha,$ with respect to the intersection form of $X,$
 on the polarization $H$. \par
For any polarized 4-manifold $(X, H)$, we can define a differential
 topological invariant \cite{w, leb-2} of $(X, H),$ by using
 Seiberg-Witten monopole equations \cite{w}. 
 We briefly recall the definition, referring to \cite{w, leb-2} for more details. Let
 ${\mathfrak{s}}$ be a spin${}^{c}$ structure of the polarized
 4-manifold $(X, H)$. Let $c_{1}({\cal L}_{\mathfrak{s}}) \in H^2(X,
 {\mathbb R})$ 
 be the first Chern class of the complex line 
bundle ${\cal L}_{\mathfrak{s}}$ associated to $\mathfrak{s}$. Suppose 
that $d_{\mathfrak{s}}:={(c^{2}_{1}({\cal L}_{\mathfrak{s}}) -2 \chi(X)
 - 3 \tau(X))}/4 =0$, which forces the virtual dimension of 
the Seiberg-Witten moduli space to be zero. Let $g$ be a $H$-adapted 
metric and assume that $c^+_{1}({\cal L}_{\mathfrak{s}}) \not=0$ with respect 
to $H={\mathcal H}^+_{g}$ is satisfied.  
Then \cite{w, leb-2}, the Seiberg-Witten invariant 
$SW_{X}(\mathfrak{s}, H)$ is defined to be the number of solutions of a 
generic perturbation of the 
Seiberg-Witten monopole equation,
modulo gauge transformation and counted with orientations. 

We can still 
define the Seiberg-Witten adapted invariant of $(X, H)$ for 
any spin${}^{c}$ structure $\mathfrak{s}$ for which $c^+_{1}({\cal L}_{\mathfrak{s}}) \not=0$ 
 and $d_{\mathfrak{s}}$ is even and positive.  In this case,  $SW_{X}(\mathfrak{s}, H)$ is defined
as the pairing $<\eta^{\frac{d_{\mathfrak{s}}}{2}},  [{\cal M}_{\mathfrak{s}}]>.$ 
Here $\eta$ is the first Chern class of the based moduli space as 
a $S^1$-bundle over the Seiberg-Witten moduli space ${\cal M}_{\mathfrak{s}}$  
and $ [{\cal M}_{\mathfrak{s}}]$ is the fundamental homology class of 
${\cal M}_{\mathfrak{s}}$. 

Hence, the Seiberg-Witten invariant 
$SW_{X}(\mathfrak{s}, H)$ of a polarized 4-manifold $(X, H)$ is 
well-defined for any spin${}^{c}$ structure $\mathfrak{s}$ with 
$c^+_{1}({\cal L}_{\mathfrak{s}}) \neq0$. Moreover, it is known \cite{mor}
that $SW_{X}(\mathfrak{s}, H)$ 
is independent of the choice of the polarization $H$ if 
$b^{+}(X) \geq 2,$ or $b^{+}(X) = 1$ and $2 \chi(X) + 3 \tau(X) > 0.$

One of the crucial properties of the Seiberg-Witten invariants 
above is that the non-triviality of the value $SW_{X}(\mathfrak{s}, H)$
 for a spin${}^{c}$ structure $\mathfrak{s}$ with 
$c^+_{1}({\cal L}_{\mathfrak{s}}) \neq0$ forces the existence of a 
non-trivial solution of the Seiberg-Witten monopole equations for 
any $H$-adapted metric. Using this, LeBrun \cite{leb-1, leb-4} 
proved the following, see also \cite{leb-2}:
\begin{thm}[\cite{leb-1, leb-4}]\label{C-bound}
Let $(X, H)$ be a polarized smooth compact oriented 4-manifold and 
let $\mathfrak{s}$ be a spin${}^{c}$ structure of $X$ and let 
$c^+_{1} \not=0$ be the orthogonal projection to $H$ with respect 
to the intersection form of $X$. Assume that $SW_{X}(\mathfrak{s}, H) \not=0$.
 Then, every $H$-adapted metric $g$ satisfies the following bounds: 
\begin{eqnarray}
{\int}_{X} s^{2}_{g} d{\mu}_{g} &\geq& 32{\pi}^{2} (c^+_{1})^{2},\notag \\
\frac{1}{4\pi^{2}}{\int}_{X} \Big(2|W^+_{g}|^{2} + \frac{s^{2}_{g}}{24} \Big) d{\mu}_{g} &\geq& \frac{2}{3}(c^+_{1})^{2}, \notag
\end{eqnarray}
with equality if and only if $g$ is K{\"{a}}hler-Einstein for a 
complex structure compatible with $\mathfrak{s}$ and has 
constant negative scalar curvature. Here, $s_{g}$ is the 
scalar curvature of $g$ and $W^+_{g}$ is the self-dual Weyl 
curvature of $g$. 
\end{thm}
We will use these bounds in Section \ref{sec-4} below. 

On the other hand, let $N$ be a closed oriented smooth 4-manifold 
with $b^{+}(N)=0$ and $k = b_{2}(N)$. By the celebrated result of 
Donaldson \cite{d}, there are classes ${\mathfrak{e}}_{1}, 
{\mathfrak{e}}_{2}, \cdots, {\mathfrak{e}}_{k} \in {H}^2(N, {\mathbb Z})$ 
descending to a basis of ${H}^2(N, {\mathbb Z})/\Tor$ with respect to 
which the intersection form is diagonal and ${\mathfrak e}_i^2=-1$ for all $i$. 
An element $\beta \in {H}^2(N, {\mathbb Z})$ is called characteristic if 
the intersection number $\beta \cdot x \equiv x \cdot x \ ( \bmod \ 2)$. 
If $\beta$ is characteristic, then $\beta \equiv w_{2}(N) \ (\bmod \ 2)$ 
and moreover there is a spin${}^c$ structure 
${\mathfrak{t}}$ on $N$ such that $c_{1}({\cal L}_{\mathfrak{t}}) = \beta$. 
Then modulo torsion, $\beta$ can  be written as 
 $ \sum^{k}_{i=1} a_{i}{\mathfrak{e}}_{i},$ where $a_i$ are integers. 
Let $x:=\sum^{k}_{i=1} x_{i}{\mathfrak{e}}_{i}$, where  $x_{i}$ are 
integers. Then we have $\beta \cdot x = -\sum^{k}_{i=1}{a}_{i}{x}_{i}$ 
and $x \cdot x = -\sum^{k}_{i=1}{x}^{2}_{i}$. This tells us that 
$\beta$ is characteristic if and only if the $a_{i}$ are odd integers, 
where $i= 1, \cdots, k$.  For example, we can obtain characteristic 
elements by taking as $a_{i} = \pm 1$. \par
The following result includes Lemma 1 of 
\cite{leb-2} as a special case.  
\begin{prop}\label{leb-lem-1}
Let $X$ be a closed oriented smooth 4-manifold with $b^{+}(X) \geq 1$ 
and $2 \chi(X) + 3 \tau(X) > 0$. Moreover, suppose that the Seiberg-Witten 
invariant of $X$ is non-trivial. Let $N$ be a closed oriented smooth 
4-manifold with $b_{1}(N)=b^{+}(N)=0$. Let $H$ be any 
polarization of a connected sum $M := X \# N$. Then there is a 
spin${}^c$ structure $\mathfrak{s}$ on $M$ such that 
$SW_{X}(\mathfrak{s}, H) \not=0$ and the self-dual part 
${c}^{+}_{1}$ of the first Chern class of the complex line bundle 
associated with $\mathfrak{s}$ satisfies
\begin{eqnarray}\label{de-bound}
({c}^{+}_{1})^2 \geq 2 \chi(X) + 3 \tau(X). 
\end{eqnarray}
\end{prop}

\begin{proof}
Notice that the Seiberg-Witten invariant of $X$ is well-defined 
and independent of the choice of the polarization under 
the assumption that $b^{+}(X) \geq 2$ or $b^{+}(X) = 1$ and 
$2 \chi(X) + 3 \tau(X) > 0$. Suppose that $\mathfrak{c}$ is the 
spin${}^{c}$ structure on $X$ with non-trivial Seiberg-Witten 
invariant. Let $\alpha :=c_{1}({\cal L}_{\mathfrak{c}}) \in H^2(X, {\mathbb Z})$ 
be the first Chern class of the complex line bundle ${\cal L}_{\mathfrak{c}}$ 
associated to $\mathfrak{c}$. Then, the non-triviality of Seiberg-Witten 
invariant forces the dimension $d_{\mathfrak{c}}$ of the Seiberg-Witten 
moduli space to be non-negative and we therefore have 
$\alpha^{2} \geq 2 \chi(X) + 3 \tau(X) > 0$. Moreover, since for any 
given polarization of $X,$ we have $ (\alpha^{+})^{2} \geq \alpha^{2},$  and we obtain 
\begin{eqnarray}\label{bound-1}
(\alpha^{+})^{2} \geq 2 \chi(X) + 3 \tau(X). 
\end{eqnarray}

Let ${\mathfrak{e}}_{1}, {\mathfrak{e}}_{2}, \cdots, {\mathfrak{e}}_{n} \in {H}^2(N, {\mathbb Z})$
  be cohomology classes descending to a basis of 
${H}^2(N, {\mathbb Z})/\Tor$ with respect to which the 
intersection form is diagonal, where $n=b_{2}(N).$ Let $H$ 
be a polarization of the connected sum $M := X \# N$. Choose new 
generators $\hat{\mathfrak{e}_{i}}={\pm} \mathfrak{e}_{i}$ for ${H}^2(N, {\mathbb Z})$ such that 
\begin{eqnarray}\label{bound-2}
\alpha^{+} \cdot (\hat{\mathfrak{e}_{i}})^{+} \geq 0
\end{eqnarray}
with respect to the polarization $H$. Then $\sum^{n}_{i=1}\hat{\mathfrak{e}_{i}}$ 
is a characteristic class and there is a spin${}^{c}$ structure $\mathfrak{t}$ on $N$ 
such that $c_{1}({\cal L}_{\mathfrak{t}}) = \sum^{n}_{i=1}\hat{\mathfrak{e}_{i}}$. 
Notice also that we have $c^{2}_{1}({\cal L}_{\mathfrak{t}}) = (\sum^{n}_{i=1}
\hat{\mathfrak{e}_{i}})^{2} = -n = -{b}_{2}(N)$.

Consider the spin${}^{c}$ structure $\mathfrak{s}:= \mathfrak{c} \# \mathfrak{t}$ 
on the connected sum $M:=X \# N.$ The first Chern class of the complex 
line bundle associated with $\mathfrak{s}$ satisfies
\begin{eqnarray*}
 c_{1}({\cal L}_{\mathfrak{s}}) = \alpha + c_{1}({\cal L}_{\mathfrak{t}}) =
 \alpha + \sum^{n}_{i=1}\hat{\mathfrak{e}_{i}}.,
\end{eqnarray*}
where we are using the same notation, $\alpha, c_{1}({\cal L}_{\mathfrak{t}}),\hat{\mathfrak{e}_{i}},$
 to denote the induced cohomology classes in $H^2(M,\ZZ).$
The gluing construction for the solutions of the Seiberg-Witten monopole 
equations on $M := X \# N,$ as in the proof of Theorem 3.1 in \cite{sung}
 (see also the proof of Proposition 2 in \cite{KMT} for $b^+\geq2$ ), 
tells us that $SW_{M}(\mathfrak{s}, H) \not=0$.

On the other hand, we obtain the following bound on $(c^{+}_{1})^2$ 
of the square of the orthogonal projection $c^{+}_{1}$ of 
$c_{1}({\cal L}_{\mathfrak{s}})$ into the polarization $H$: 
\begin{eqnarray*}
(c^{+}_{1})^2 &=& \Big( \alpha^{+} + \sum^{n}_{i=1} (\hat{\mathfrak{e}_{i}})^{+}
 \Big)^{2} = ( \alpha^{+})^{2} + 2 \sum^{n}_{i=1} (\alpha^{+} \cdot (\hat{\mathfrak{e}_{i}})^{+})
 + \sum^{n}_{i=1} ((\hat{\mathfrak{e}_{i}})^{+})^{2} \\
             &\geq & ( \alpha^{+})^{2}, 
\end{eqnarray*}
where we used (\ref{bound-2}). This bound and (\ref{bound-1}) implies the desired bound (\ref{de-bound}). 
\end{proof}

\begin{rmk}
Under the assumptions of Proposition \ref{leb-lem-1}, suppose that $X$ satisfies $b^{+}(X) = 1$ and $2 \chi(X) + 3 \tau(X) > 0$. Then the connected sum $M := X \# N$ satisfies $b^{+}(M) = 1$ and $2 \chi(M) + 3 \tau(M) = 2 \chi(X) + 3 \tau(X)-b_{2}(N)$.  In general, if $2 \chi(M) + 3 \tau(M) <0$ the Seiberg-Witten invariant depends of the choice of polarization. However, the above proof implies that if $M:= X \# N,$ with $X,N$ as in Proposition \ref{leb-lem-1}, any polarization on $M$ satisfies the following bound $(c^{+}_{1})^2\geq ( \alpha^{+})^{2}\geq 2 \chi(X) + 3 \tau(X) > 0$. This tells us that the value of the Seiberg-Witten invariant is independent of the choice of polarization.
\end{rmk}

\section{Non-singular solutions and a bound of the Ricci curvature}\label{sec-3}

The main result of this section is Proposition \ref{key-prop} below. We shall use it in the next section to find an obstruction to the existence of non-singular solutions of the normalized Ricci flow. \par 
We start by recalling the following result on the trace free part of the Ricci curvature of the long time solution of the normalized Ricci flow. 
\begin{lem}[\cite{fz-1, ism}]\label{FZZ-prop}
Let $X$ be a closed oriented Riemannian $n$-manifold and assume that there is a long time solution $\{g(t)\}$, $t \in [0, \infty)$, to the normalized Ricci flow. Assume moreover that the solution satisfies $|{s}_{g(t)}| \leq C$ and 
\begin{eqnarray}\label{mini-bound}
\hat{s}_{g(t)}:=\min_{x \in X}{s}_{g(t)}(x) \leq -c <0,
\end{eqnarray}
where the constants $C$ and $c$ are independent of both $x \in X$ and time $t \in [0, \infty)$. Then, the trace-free part $\stackrel{\circ}{r}_{g(t)}$ of the Ricci curvature satisfies 
\begin{eqnarray*}
{\int}^{\infty}_{0} {\int}_{X} |\stackrel{\circ}{r}_{g(t)}|^2 d{\mu}_{g(t)}dt < \infty. 
\end{eqnarray*}
In particular, as $m \rightarrow \infty$
\begin{eqnarray*}
{\int}^{m+1}_{m} {\int}_{X} |\stackrel{\circ}{r}_{g(t)}|^2 d{\mu}_{g(t)}dt \longrightarrow 0. 
\end{eqnarray*}
\end{lem}

On the other hand, there is a natural diffeomorphism invariant arising from a variational problem of the total scalar curvature of Riemannian metrics on any given closed oriented Riemannian manifold $X$ of dimension $n\geq 3$. As was conjectured by Yamabe \cite{yam}, and later proved by Trudinger, Aubin, and Schoen \cite{aubyam, rick, trud}, every conformal class on any smooth compact manifold contains a Riemannian metric of constant scalar curvature.
To be more precise, for any conformal class $[g]=\{ vg ~|~v: X\to {\Bbb R}^+\}$, we can consider an associated number $Y_{[g]}$ which is called the { Yamabe constant} of the conformal class $[g]$ and is defined by 
\begin{eqnarray*}
Y_{[g]} = \inf_{h \in [g]}  \frac{\int_X 
s_{{h}}~d\mu_{{h}}}{\left(\int_X 
d\mu_{{h}}\right)^{\frac{n-2}{n}}}, 
\end{eqnarray*}
where $d\mu_{{h}}$ is the volume form with respect to the metric $h.$ It is known \cite{aubyam, rick, trud}  that this number is realized as the constant scalar curvature of some  metric in the conformal class $[g].$ Then, Kobayashi \cite{kob} and Schoen \cite{sch} independently introduced the following invariant of 
$X$:
\begin{eqnarray*}
{\mathcal Y}(X) = \sup_{[g] \in \mathcal{C}}Y_{[g]}, 
\end{eqnarray*}
where $\mathcal{C}$ is the set of all conformal classes on $X$. This is now known as the {Yamabe invariant} of $X$. We have the following bound:
\begin{lem}[\cite{ism}]\label{yama-b}
Let $X$ be a closed oriented Riemannian manifold of dimension $n \geq 3$ and assume that the Yamabe invariant of $X$ is negative, i.e., ${\mathcal Y}(X)<0$. If there is a solution $\{g(t)\}$, $t \in [0, T)$, to the normalized Ricci flow, then the solution satisfies the bound (\ref{mini-bound}). More precisely, the following is satisfied:
\begin{eqnarray*}
\hat{s}_{g(t)}:=\min_{x \in X}{s}_{g(t)}(x) \leq \frac{{\mathcal Y}(X)}{(vol_{g(0)})^{2/n}} < 0. 
\end{eqnarray*}
\end{lem}
We recall now the following definition: 
\begin{defn}[\cite{ism}]\label{bs}
A maximal solution $\{g(t)\}$, $t \in [0, T)$, to the normalized Ricci flow  on $X$ is called quasi-non-singular if $T=\infty$ and the scalar curvature $s_{g(t)}$ of $g(t)$ satisfies 
\begin{eqnarray*}
\sup_{X \times [0, \infty)}|{s}_{g(t)}| < \infty. 
\end{eqnarray*}
\end{defn}
Notice that any non-singular solution is quasi-non-singular.
\begin{prop}\label{key-prop}
Let $X$ be a closed oriented smooth 4-manifold with $b^{+}(X) \geq 1$ and $2 \chi(X) + 3\tau(X) > 0$. Assume that $X$ has a non-trivial Seiberg-Witten invariant. Let $N$ be a closed oriented smooth 4-manifold with $b_{1}(N)=b^{+}(N)=0$. If there is a quasi-non-singular solution to the normalized Ricci flow  on the connected sum $M:=X \# N$, then the trace-free part $\stackrel{\circ}{r}_{g(t)}$ of the Ricci curvature satisfies 
\begin{eqnarray}\label{fzz-ricci-0112}
{\int}^{\infty}_{0} {\int}_{X} |\stackrel{\circ}{r}_{g(t)}|^2 d{\mu}_{g(t)}dt < \infty. 
\end{eqnarray}
In particular, 
\begin{eqnarray}\label{fzz-ricci-011}
{\int}^{m+1}_{m} {\int}_{X} |\stackrel{\circ}{r}_{g(t)}|^2 d{\mu}_{g(t)}dt \longrightarrow 0
\end{eqnarray}
holds when $m \rightarrow +\infty$. 
\end{prop}

\begin{proof}
First of all, notice that the connected sum $M$ has non-trivial Seiberg-Witten invariant with respect to any polarization by Proposition \ref{leb-lem-1}. By Witten's vanishing theorem \cite{w}, this implies that $M$ cannot admit any metric of positive scalar curvature. On the other hand, it is known \cite{leb-3} that the Yamabe invariant of any closed $n$-manifold $Z$ which cannot admit metrics of positive scalar curvature is given by
\begin{eqnarray}\label{yama-1}
{\mathcal Y}(Z) = - \Big(\inf_{g \in {\cal R}_{Z}}{\int}_{Z}|s_{g}|^{{n}/{2}} d{\mu}_{g} \Big)^{{2}/{n}}, 
\end{eqnarray}
where $ {\cal R}_{Z}$ is the set of all Riemannian metrics on $Z$. Combining the first inequality in Theorem \ref{C-bound} and the inequality in Proposition \ref{leb-lem-1} we get the following bound: 
\begin{eqnarray}\label{yama-2}
 {\int}_{M} s^{2}_{g} d{\mu}_{g} \geq 32{\pi}^{2} \Big( 2 \chi(X) + 3\tau(X) \Big). 
\end{eqnarray}
Note that this bound holds for any $H$-adapted metric on $M$, where $H$ is any polarization of $M.$ In particular, it holds for any metric $g$. 

Moreover, (\ref{yama-1}) and (\ref{yama-2}) implies 
\begin{eqnarray}\label{yama-3}
{\mathcal Y}(M) \leq -4{\pi}\sqrt{2(2 \chi(X) + 3\tau(X))} < 0. 
\end{eqnarray}
This bound and Lemma \ref{yama-b} tell us that any solution to the normalized Ricci flow on $M$ satisfies
\begin{eqnarray*}
\hat{s}_{g(t)}:=\min_{x \in M}{s}_{g(t)}(x) \leq {-4{\pi}} \sqrt{ \frac{2(2 \chi(X) + 3\tau(X))}{vol_{g(0)}}} < 0. 
\end{eqnarray*}

But, this last inequality combined with Lemma \ref{FZZ-prop} shows that any quasi-non-singular solution to the normalized Ricci flow on the connected sum $M$  must satisfy the desired bound (\ref{fzz-ricci-0112}). 
\end{proof}

\section{Proof of Theorem \ref{main-A}}\label{sec-4}

We are now in position to prove Theorem \ref{main-A}, which is a special case of the following result: 
\begin{thm}\label{main-C}
Let $X$ be a closed oriented smooth 4-manifold with $b^+(X) \geq 1$ and $2 \chi(X) + 3\tau(X) > 0$. Assume that $X$ has a non-trivial Seiberg-Witten invariant. Let $N$ be a closed oriented smooth 4-manifold with $b_{1}(N)=b^{+}(N)=0$. Then, there do not exist quasi-non-singular solutions to the normalized Ricci flow on $M:=X \# N$ if 
\begin{eqnarray}\label{ob-N-Ricci-2}
b_{2}(N) > \frac{1}{3} \Big( 2\chi(X) + 3 \tau(X) \Big). 
\end{eqnarray}
In particular, there is no non-singular solutions to the normalized Ricci flow on $M$. \par 
Moreover, if we assume that $N$ is not an integral homology 4-sphere whose fundamental group has no non-trivial finite quotient, then there do not exist quasi-non-singular solutions to the normalized Ricci flow on $M$ if 
\begin{eqnarray}\label{ob-N-Ricci-3}
b_{2}(N) \geq  \frac{1}{3} \Big( 2\chi(X) + 3 \tau(X) \Big). 
\end{eqnarray}
In particular, there are no non-singular solutions to the normalized Ricci flow on $M$. 
\end{thm}
\begin{proof}
Suppose that there would be a quasi-non-singular solution $\{g(t) \}$ to the normalized Ricci flow on the connected sum $M:=X \# N$. Then the second inequality in Theorem \ref{C-bound} tells us that, for any time $t$,  $g(t)$ must satisfy 
\begin{eqnarray*}
\frac{1}{4\pi^{2}}{\int}_{X} \Big(2|W^+_{g(t)}|^{2} + \frac{s^{2}_{g(t)}}{24} \Big) d{\mu}_{g(t)} &\geq& \frac{2}{3}(c^+_{1})^{2}. 
\end{eqnarray*}
for any spin${}^{c}$ structure $\mathfrak{s}$ with $SW_{M}(\mathfrak{s}, H) \not=0$, where $H:={\cal H}^{+}_{g(t)}$. However, Proposition \ref{leb-lem-1} now asserts that the connected sum $M:=X \# N$ has a spin${}^{c}$ structure with $({c}^{+}_{1})^2 \geq 2 \chi(X) + 3 \tau(X)$. We therefore conclude that 
\begin{eqnarray}\label{mono-1}
\frac{1}{4{\pi}^2}{\int}_{M}\Big(2|W^{+}_{g(t)}|^2+\frac{{s}^2_{g(t)}}{24}\Big) d{\mu}_{g(t)} \geq  \frac{2}{3}\Big( 2 \chi(X) + 3\tau(X)\Big). 
\end{eqnarray}
Moreover, the equality holds 
if and only if the metric $g(t)$ is a K{\"{a}}hler-Einstein metric with negative scalar curvature \cite{leb-4}. 

On the other hand, we have the following Gauss-Bonnet like formula:
\begin{eqnarray*}\label{GB}
2\chi(M) + 3\tau(M) = \frac{1}{4{\pi}^2}{\int}_{M}\Big(2|W^{+}_{g(t)}|^2+\frac{{s}^2_{g(t)}}{24}-\frac{|{r}^{\circ}_{g(t)}|^2}{2} \Big) d{\mu}_{g(t)}. 
\end{eqnarray*}
In particular, we obtain
\begin{eqnarray*}
2\chi(M) + 3\tau(M) &=& {\int}^{m+1}_{m} \Big(2\chi(M) + 3\tau(M) \Big)dt \\
&=& \frac{1}{4{\pi}^2}{\int}^{m+1}_{m} {\int}_{M}\Big(2|W^{+}_{g(t)}|^2+\frac{{s}^2_{g(t)}}{24}-\frac{|{r}^{\circ}_{g(t)}|^2}{2} \Big) d{\mu}_{g(t)}dt.
\end{eqnarray*}
Since Proposition \ref{key-prop} tells us that a quasi-non-singular solution $\{g(t)\}$ must satisfy (\ref{fzz-ricci-011}), by taking ${m \longrightarrow \infty}$ in the above inequality, we obtain 
\begin{eqnarray}\label{key-11}
2\chi(M) + 3\tau(M) = \lim_{m \longrightarrow \infty}\frac{1}{4{\pi}^2}{\int}^{m+1}_{m} {\int}_{M}\Big(2|W^{+}_{g(t)}|^2+\frac{{s}^2_{g(t)}}{24} \Big) d{\mu}_{g(t)}dt.
\end{eqnarray}
Moreover, by the inequality (\ref{mono-1}), we get 
\begin{eqnarray*}
\frac{1}{4{\pi}^2}{\int}^{m+1}_{m} {\int}_{M}\Big(2|W^{+}_{g(t)}|^2+\frac{{s}^2_{g(t)}}{24}\Big) d{\mu}_{g(t)}dt &\geq& \frac{2}{3}{\int}^{m+1}_{m} \Big( 2 \chi(X) + 3\tau(X)\Big) dt \\
&=& \frac{2}{3}\Big( 2 \chi(X) + 3\tau(X)\Big).
\end{eqnarray*}
This bound and (\ref{key-11}) tell us that the following holds: 
\begin{eqnarray*}
2\chi(M) + 3\tau(M) \geq \frac{2}{3}\Big( 2 \chi(X) + 3\tau(X)\Big).  
\end{eqnarray*}
Since we have $2\chi(M) + 3\tau(M) = 2\chi(X) + 3 \tau(X) - {b}_{2}(N)$, we get 
\begin{eqnarray*}
 2\chi(X) + 3 \tau(X) - {b}_{2}(N)  \geq  \frac{2}{3}\Big( 2 \chi(X) + 3 \tau(X) \Big)
\end{eqnarray*}
which is equivalent to
\begin{eqnarray*}
{b}_{2}(N) \leq \frac{1}{3}\Big( 2\chi(X) + 3 \tau(X) \Big). 
\end{eqnarray*}

By contraposition, we conclude that under (\ref{ob-N-Ricci-2}), there do not exist quasi-non-singular solutions on $M$. In particular, there is no non-singular solutions to the normalized Ricci flow. \par
Suppose moreover that $N$ is not an integral homology 4-sphere whose fundamental group has no non-trivial finite quotient. Then we observe the equality cannot occur in (\ref{mono-1}). Notice that the equality case in (\ref{mono-1}) forces that $g(t)$ must be a K{\"{a}}hler-Einstein metric with negative scalar curvature. In particular, this forces that the connected sum $M:=X \# N$ is a minimal K{\"{a}}hler surface \cite{a, yau}. On the other hand, Theorem 5.4 in \cite{k} ($b^{+}(M) \geq 2$ case) and Theorem 2 in \cite {h-k} ($b^{+}(M) = 1$ case) tell us that if a minimal K{\"{a}}hler surface admits a connected sum decomposition $X \# N$, then $N$ must be an integral homology 4-sphere whose fundamental group has no non-trivial finite quotient. Therefore, we conclude that the equality cannot occur in (\ref{mono-1}). Hence we have the following strict inequality which holds for $g(t)$ on $M$: 
\begin{eqnarray}\label{mono-2}
\frac{1}{4{\pi}^2}{\int}_{M}\Big(2|W^{+}_{g(t)}|^2+\frac{{s}^2_{g(t)}}{24}\Big) d{\mu}_{g(t)} > \frac{2}{3}\Big( 2 \chi(X) + 3\tau(X)\Big). 
\end{eqnarray}
This strict inequality and the strategy of the proof of bound (\ref{ob-N-Ricci-2}) above enable us to conclude us that under (\ref{ob-N-Ricci-3}), there does not exist quasi-non-singular solutions on $M$. 
\end{proof}
\begin{rmk}
Notice that any Einstein metric is a fixed point of the normalized Ricci flow. In particular, we can obtain a non-singular solution of the normalized Ricci flow by taking an Einstein metric as a initial metric of the normalized Ricci flow. Therefore, Theorem \ref{main-C} recovers 
the obstruction to the existence of Einstein metrics, due to LeBrun \cite{leb-4}. 
\end{rmk}

\section{Proof of Theorem \ref{main-B}}\label{sec-5}

On Del Pezzo surfaces, we have a complete solution for the existence of K{\"a}hler-Einstein metrics, due to Tian:
\begin{thm}[\cite{tian}]\label{positive}
A compact complex surface $X$ admits a K{\"{a}}hler-Einstein metric with positive scalar curvature if and only if its anti-canonical line bundle is ample and its Lie algebra of holomorphic vector fields is reductive. 
\end{thm}
By this theorem, the 
complex surfaces ${\mathbb C}{P}^2 \# {\ell} \overline{{\mathbb C}{P}^2}$, where $3 \leq \ell \leq 8$, admit K{\"{a}}hler-Einstein metrics with positive scalar curvature. 

In more generality, the normalized Ricci flow on compact K{\"a}hler manifolds was studied by Cao  \cite{c, c-c} . We are going to use the following result:

\begin{thm}[\cite{c, c-c}]\label{cao-K}
Let $M$ be a compact K{\"{a}}hler manifold with definite first Chern class ${c}_{1}(M)$. If ${c}_{1}(M)=0$, then for any initial K{\"{a}}hler metric $g_{0}$, the solution to the normalized Ricci flow exists for all time and converges to a Ricci-flat metric as $t \rightarrow \infty$. If ${c}_{1}(M) < 0$ and the initial metric $g_0$ is chosen to represent minus the first Chern class, the solution to the normalized Ricci flow exists for all time and converges to an Einstein metric of negative scalar curvature as $t \rightarrow \infty$. If ${c}_{1}(M) > 0$ and the initial metric $g_0$ is chosen to represent the first Chern class, then the solution to the normalized Ricci flow exists for all time. 
\end{thm}
Notice that, in case when ${c}_{1}(M) = 0$ or ${c}_{1}(M) < 0$, the solution is actually non-singular \cite{c-c} in the sense of Definition \ref{non-sin}. On the other hand, the non-triviality of the Seiberg-Witten invariants in the case when  ${c}_{1}(M) < 0$ tells us that the scalar curvature $s_{g_{0}}$ of the initial metric $g_{0}$ cannot have $s_{g_{0}} \geq 0$. \par

On complex surfaces of general type, existence results for simply connected manifolds with $b^+=1$ and ample canonical line bundle have only been recently found. One of the first examples is a deformation of the Barlow surface for which the ampleness of the canonical bundle is due to Catanese and LeBrun \cite{CL}:
\begin{thm}[\cite{CL}]\label{ample0}
There exist simply connected complex surfaces of general type, with $b^+=1$, $c_{1}^{2}=1$ and ample canonical bundle.
\end{thm}
In \cite{R-S}, the second and third authors of this paper proved the following result:
\begin{thm}[\cite{R-S}]\label{ample}
There exist simply connected complex surfaces of general type, with $b^+=1$, $c_{1}^{2}=2$ or $3$ and ample canonical bundle. 
\end{thm}
Using the same method as in \cite{R-S} Park and his collaborators \cite{pps}, showed that the exotic manifold that they constructed also admits an ample canonical line bundle:
\begin{thm}[\cite{pps}]\label{ample'}
There exist simply connected complex surfaces of general type, with $b^+=1$, $c_{1}^{2}=4$ and ample canonical bundle. 
\end{thm}

We are now ready to prove Theorem \ref{main-B}:\par
\begin{proof}
If we want to construct a smooth structure on the manifold $X$ which has positive Yamabe invariant and admits solutions for the normalized Ricci flow, then we can just consider 
the canonical differential and complex structures of the complex projective plane blown-up at $l$ points, where $3\leq l\leq 8.$ 
The existence of  an Einstein metric metric is given by Theorem \ref{positive}.
Hence, on these 
Del Pezzo surfaces, there are non-singular solutions (fixed points) of the normalized Ricci flow  by taking the K{\"{a}}hler-Einstein metrics with positive scalar curvature as initial metrics. Since the scalar curvature of their metrics is positive, we notice that according to Lemma 1.5 in \cite{kob}, the Yamabe invariant of these manifolds must be positive. \par
In the second case of the theorem, we are going to consider the smooth structures associated to the complex structures of general type found in Theorems \ref{ample0}, \ref{ample}, \ref{ample'}. On these manifolds, Cao's Theorem \ref{cao-K} tells us that solutions to the normalized Ricci flow exist if we start with a  K{\"a}hler metric whose K{\"a}hler form is in the cohomology class of the canonical line bundle. Moreover, for surfaces of general type  the Yamabe invariant \cite{leb-3} is strictly negative. \par
For the third part of the proof, we are going to use  the constructions in Section 3 of \cite{R-S}. In  the proof of Proposition 3.2 in \cite{R-S}, it is showed that for any $M:=\bcp^2\#\ell \cpb, 5\leq \ell\leq 8,$ there exist infinitely many manifolds of the form $M_i:=X_i\#(\ell-3)\cpb, i\in\NN,$ which are homeomorphic, non-diffeomorphic to $M$ and to each other. Here, $X_i$ are the 4-manifolds  homeomorphic to $\bcp^2\#3\cpb$ constructed in \cite{ABP} and which have non-trivial Seiberg-Witten invariant. Then the manifolds $M_i$ have non-trivial Seiberg-Witten invariant, and as $c_1^2(M_i)=9-l>0,$ the Yamabe invariant ${\mathcal Y }(M_i)$ is strictly negative \cite{leb-3}.  The manifolds $X_i$ have $c_1^2(X_i)=6,$ non-trivial Seiberg-Witten invariant by construction and of course $(\ell-3)\geq\frac 13c_1^2(X_i)=2,$ as $5\leq \ell\leq 8.$ Hence, Theorem \ref{main-A} tells us that there are no solutions to the normalized
  Ricci flow on any $M_i$ for any initial metric.
\end{proof}


\vfill

{\footnotesize 
\noindent
{Masashi Ishida,} \\
{Department of Mathematics,  
Sophia University, \\ 7-1 Kioi-Cho, Chiyoda-Ku, 
 Tokyo 102-8554, Japan }\\
{\sc e-mail}: ishida@mm.sophia.ac.jp}

\vspace{0.1cm}

{\footnotesize 
\noindent
{Rare{\c s} R{\u a}sdeaconu}, \\
{Einstein Institute of Mathematics\\
Hebrew University of Jerusalem\\
Givat Ram, Jerusalem, 91904, Israel}\\
{\sc e-mail}: rares@math.huji.ac.il}

\vspace{0.1cm}

{\footnotesize 
\noindent
{Ioana {\c S}uvaina}, \\
{Courant Institute of Mathematical Sciences, \\
251 Mercer St. New York, 10012, USA}\\
{\sc e-mail}: ioana@cims.nyu.edu}

\end{document}